\documentclass[12pt,reqno]{amsart}
\usepackage{enumerate, latexsym, amsmath, amsfonts, amssymb, amsthm, color}
\def\Z{\Bbb Z}
\def\N{\Bbb N}

\def\l{\left}
\def\r{\right}
\def\bg{\bigg}
\def\({\bg(}
\def\){\bg)}
\def\t{\text}
\def\f{\frac}
\def\mo{{\rm{mod}\ }}

\def\sgn{{\rm sgn}}
\def\cs{\ldots}

\def\ls{\leqslant}
\def\gs{\geqslant}
\def\se {\subseteq}
\def\sm{\setminus}
\def\bi{\binom}

\def\eq{\equiv}

\def\da{\delta}

\def\colon{{:}\;}
\def\Proof{\noindent{\it Proof}}
\def\Ack{\medskip\noindent {\bf Acknowledgments}}
\theoremstyle{plain}
\newtheorem{theorem}{Theorem}

\newtheorem{lemma}{Lemma}
\newtheorem{corollary}{Corollary}
\newtheorem{conjecture}{Conjecture}
\theoremstyle{definition}

\theoremstyle{remark}
\newtheorem{remark}{Remark}

 \vspace{4mm}

\begin{document}
 \baselineskip=17pt
\hbox{J. Combin. Theory Ser. A 119(2012), no.\,2, 364--381.}
\medskip

\title
[Linear extension of the Erd\H os-Heilbronn conjecture] {Linear
extension of the Erd\H os-Heilbronn conjecture}

\author
[Zhi-Wei Sun and Li-Lu Zhao] {Zhi-Wei Sun* and Li-Lu Zhao}

\thanks{*Supported by the National Natural Science Foundation (grant 10871087 and 11171140)
 of China}

\address {(Zhi-Wei Sun) Department of Mathematics, Nanjing
University, Nanjing 210093, People's Republic of China}
\email{zwsun@nju.edu.cn}

\address{(Li-Lu Zhao) Department of Mathematics
\\The University of Hong Kong
\\Pokfulam Road, Hong Kong}
\email{zhaolilu@gmail.com}


 \begin{abstract}  The famous Erd\H os-Heilbronn conjecture plays an
important role in the development of additive combinatorial number
theory. In 2007 Z. W. Sun made the following further conjecture
(which is the linear extension of the Erd\H os-Heilbronn
conjecture): For any finite subset $A$ of a field $F$ and nonzero
elements $a_1,\ldots,a_n$ of $F$, we have
$$\aligned&|\{a_1x_1+\cdots+a_nx_n:\ x_1,\ldots,x_n\in A,\ \t{and}\ x_i\not= x_j\ \t{if}\ i\not=j\}|
\\&\qquad\qquad\gs\min\l\{p(F)-\da,\ n(|A|-n)+1\r\},
\endaligned$$
where the additive order $p(F)$ of the multiplicative identity of
$F$ is different from $n+1$, and $\da\in\{0,1\}$ takes the value 1 if
and only if $n=2$ and $a_1+a_2=0$. In this paper we prove this
conjecture of Sun when $p(F)\gs n(3n-5)/2$. We also obtain a sharp
lower bound for the cardinality of the restricted sumset
$$\{x_1+\cdots+x_n:\ x_1\in A_1,\ldots,x_n\in A_n,\ \t{and}\ P(x_1,\ldots,x_n)\not=0\},$$
where $A_1,\ldots,A_n$ are finite subsets of a field $F$ and
$P(x_1,\ldots,x_n)$ is a general polynomial over $F$.
\end{abstract}

\maketitle

\section{Introduction}
\setcounter{lemma}{0}
\setcounter{theorem}{0}
\setcounter{corollary}{0}
\setcounter{remark}{0}
\setcounter{equation}{0}
\setcounter{conjecture}{0}

A basic objective in the active field of additive combinatorial
number theory is the sumset of finite subsets $A_1,\ldots,A_n$ of a
field $F$ given by
$$A_1+\cdots+A_n=\{x_1+\cdots+x_n:\ x_1\in A_1,\ldots,x_n\in A_n\}.$$
(See, e.g., \cite{N96} and \cite{TV}.) The well-known Cauchy-Davenport theorem
asserts that
$$|A_1+\cdots+A_n|\gs\min\{p(F),\ |A_1|+\cdots+|A_n|-n+1\},$$
where $p(F)$ is the additive order of the multiplicative identity of
$F$ (which is the characteristic of $F$ if $F$ is of a prime
characteristic, and the positive infinity if $F$ is of
characteristic zero). When $n=2$ and $F=\Z/p\Z$ with $p$ a prime,
this gives the original form of the Cauchy-Davenport theorem.

In 1964 P. Erd\H os and H. Heilbronn \cite{EH} conjectured that if $p$ is
a prime and $A$ is a subset of $\Z/p\Z$ then
$$|\{x+y:\ x,y\in A\ \t{and}\ x\not=y\}|\gs\min\{p,\ 2|A|-3\}.$$
This challenging conjecture was finally solved by J. A. Dias da
Silva and Y. O. Hamidoune \cite{DH} in 1994 who employed exterior
algebras to show that for any subset $A$ of a field $F$ we have
$$|\{x_1+\cdots+x_n:\, x_i\in A,\ \ x_i\not=x_j\ \t{if}\ i\not=j\}|
\gs \min\{p(F),\, n|A|-n^2+1\}.$$ Recently P. Balister and J. P.
Wheeler \cite{BW} extended the Erd\H os-Heilbronn conjecture to any
finite group.

In 1995-1996 N. Alon, M. B. Nathanson and I. Z. Ruzsa \cite{ANR1}\cite{ANR2}
used the so-called polynomial method rooted in \cite{AT} to prove that if
$A_1,\ldots,A_n$ are finite subsets of a field $F$ with
$0<|A_1|<\cdots<|A_n|$ then
$$|\{x_1+\cdots+x_n:\,  x_i\in A_i,\ x_i\not=x_j\ \t{if}\ i\not=j\}|
\gs\min\bg\{p(F),\sum_{i=1}^n(|A_i|-i)+1\bg\}.$$ The polynomial
method was further refined by Alon \cite{A99} in 1999, who presented the
following useful principle.

\medskip
\noindent {\bf Combinatorial Nullstellensatz} {\rm (Alon \cite{A99})}.
{\it Let $A_1,\ldots,A_n$ be finite subsets of a field $F$ with
$|A_i|>k_i$ for all $i=1,\ldots,n$ where
$k_1,\ldots,k_n\in\N=\{0,1,2,\ldots\}$. Suppose that
$P(x_1,\ldots,x_n)$ is a polynomial over $F$
 with  $[x_1^{k_1}\cdots x_n^{k_n}]P(x_1,\ldots,x_n)$ $($the coefficient of the monomial
$x_1^{k_1}\cdots x_n^{k_n}$ in the polynomial $P(x_1,\ldots,x_n))$
nonzero and $k_1+\cdots+k_n=\deg P$. Then there are $x_1\in
A_1,\ldots,x_n\in A_n$ such that $P(x_1,\ldots,x_n)\not=0$.}
\medskip

The Combinatorial Nullstellensatz has been applied to investigate
some sumsets with polynomial restrictions by various authors, see
\cite{DKSS,HS,LS,PS1,S03,SY,K05,S08b}.

Throughout this paper, for a predicate $P$ we let
$$[\![P]\!]=\begin{cases}1&\t{if}\ P\ \t{holds},\\0&\t{otherwise}.\end{cases}$$
For $a,b\in\Z$ we define $[a,b]=\{m\in\Z:\ a\ls m\ls b\}$. For a
field $F$ we let $F^*$ be the multiplicative group of all nonzero
elements of $F$. As usual the symmetric group on $\{1,\ldots,n\}$ is
denoted by $S_n$. For $\sigma\in S_n$ we use $\sgn(\sigma)$ to stand
for the sign of the permutation $\sigma$. We also set $(x)_0=1$ and
$(x)_n=\prod_{j=0}^{n-1}(x-j)$ for $n=1,2,3,\ldots$.

Recently Z. W. Sun  made the following conjecture (cf. \cite{S08a}) which can be
viewed as the linear extension of the Erd\H os-Heilbronn conjecture.

\begin{conjecture} [Sun] \label{Conj1.1} Let $A$ be a finite
subset of a field $F$ and let $a_1,\ldots,a_n\in F^*=F\sm\{0\}$.
Provided  $p(F)\not=n+1$ we have
\begin{equation}\label{1.1}\begin{aligned}&|\{a_1x_1+\cdots+a_nx_n:\ x_1,\ldots,x_n\in A,\ \t{and}\ x_i\not= x_j\ \t{if}\ i\not=j\}|
\\&\qquad\gs\min\l\{p(F)-[\![n=2\ \&\ a_1=-a_2]\!],\ n(|A|-n)+1\r\}.
\end{aligned}\end{equation}
\end{conjecture}

 \medskip

\noindent{\it Example} 1.1. Let $p$ be an odd prime and let $k$ be a positive integer relatively prime to $p-1$.
As $k\not\eq0\ (\mo\ p-1)$, we have $\sum_{x\in F_p}x^k=0$ where $F_p=\Z/p\Z$.
For any distinct $x,y\in F_p$ we cannot have $x^k=y^k$ since $ku+(p-1)v=1$ for some $u,v\in\Z$.
Thus
$$|\{x_1^k+\cdots+x_{p-2}^k+2x_{p-1}^k:\ x_1,\ldots,x_{p-1}\in F_p\ \t{are
distinct}\}|=|F_p^*|=p-1.$$
In the case $k=1$ and $p\in\{5,7\}$, this was noted by Mr. Wen-Long
Zhang (at Nanjing University) in Feb. 2011 via computation under the guidance of the first author.
\medskip

All known proofs of the Erd\H os-Heilbronn conjecture (including the
recent one given by S. Guo and Sun \cite{GS} based on Tao's harmonic
analysis method) cannot be modified easily to confirm the above
conjecture. New ideas are needed!

Concerning Conjecture 1.1 we are able to establish the following
result.

\begin{theorem} \label{Th1.1} Let $A$ be a finite subset of a field $F$
and let $a_1,\ldots,a_n\in F^*$. Then $(1.1)$ holds if $p(F)\gs
n(3n-5)/2$.
\end{theorem}

We obtain Theorem 1.1 by combining our next two theorems.

\begin{theorem} \label{Th1.2} Let $n$ be a positive integer, and let $F$ be a
field with $p(F)\gs(n-1)^2$. Let $a_1,\ldots,a_n\in F^*$, and
suppose that $A_i\se F$ and $|A_i|\gs 2n-2$ for $i=1,\ldots,n$.
Then, for the set
\begin{equation}\label{1.2}C=\{a_1x_1+\cdots+a_nx_n:\ x_1\in A_1,\ldots,x_n\in A_n,\ \t{and}\ x_i\not=x_j\ \t{if}\ i\not=j\}
\end{equation}
we have
\begin{equation}\label{1.3}|C|\gs\min\{p(F)-[\![n=2\ \&\ a_1+a_2=0]\!],\,|A_1|+\cdots+|A_n|-n^2+1\}.
\end{equation}
\end{theorem}

Theorem \ref{Th1.2} has the following consequence.

\begin{corollary} \label{Cor1.1} Let $p>7$ be a prime and let $A\se
F_p=\Z/p\Z$ with $|A|\gs\sqrt{4p-7}$. Let $n=\lfloor|A|/2\rfloor$
and $a_1,\ldots,a_n\in F_p^*$. Then every element of $F_p$ can be
written in the linear form $a_1x_1+\cdots+a_nx_n$ with
$x_1,\ldots,x_n\in A$ distinct.
\end{corollary}

\begin{remark} In the case $a_1=\cdots=a_n=1$, Corollary \ref{Cor1.1} is a
refinement of a conjecture of Erd\H os proved by da Silva and
Hamidoune \cite{DH} via exterior algebras.
\end{remark}

By Theorem \ref{Th1.1}, Conjecture \ref{Conj1.1} is valid for $n=2$. Now we explain
why Conjecture \ref{Conj1.1} holds in the case $n=3$. Let $A$ be a finite
subset of a field $F$ and let $a_1,a_2,a_3\in F^*$. Clearly (\ref{1.1})
holds if $|A|\ls n$. Below we assume $|A|>n=3$.
By Theorem \ref{Th1.1}, (\ref{1.1}) with $n=3$ holds if
$p(F)\gs 3(3\times 3-5)/2=6$.  When $p(F)=5$, we have (\ref{1.1}) by Theorem \ref{Th1.2}.
If $p(F)=2$
and $c_1,c_2,c_3,c_4$ are four distinct elements of $A$, then
\begin{align*}&|\{a_1x_1+a_2x_2+a_3x_3:\ x_1,x_2,x_3\in A\ \t{and}\
x_1,x_2,x_3\ \t{are distinct}\}|
\\&\qquad\gs|\{a_1c_1+a_2c_2+a_3c_3,\,a_1c_1+a_2c_2+a_3c_4\}|
\\&\qquad=2=\min\{p(F),\,3(|A|-3)+1\}.
\end{align*}
In the case $p(F)=3$,  for some $1\ls s<t\ls 3$ we have
$a_s+a_t\not=0$, hence for any $c\in A$ we have
\begin{align*}&|\{a_1x_1+a_2x_2+a_3x_3:\ x_1,x_2,x_3\in A\ \t{and}\
x_1,x_2,x_3\ \t{are distinct}\}|
\\\gs&|\{a_sx_s+a_tx_t:\ x_s,x_t\in A\sm\{c\}\ \t{and}\ x_s\not=x_t\}|
\\\gs&\min\{p(F),\,2(|A\sm\{c\}|-2)+1\}\ \ (\t{by Theorem 1.1 with}\ n=2)
\\=&3=\min\{p(F),\,3(|A|-3)+1\}.
\end{align*}

In this paper we also apply the Combinatorial Nullstellensatz {\it
twice} to deduce the following result on sumsets with {\it general}
polynomial restrictions.

\begin{theorem}\label{Th1.3} Let $P(x_1,\ldots,x_n)$ be a polynomial
over a field $F$. Suppose that $k_1,\ldots,k_n$ are nonnegative
integers with $k_1+\cdots+k_n=\deg P$ and $[x_1^{k_1}\cdots
x_n^{k_n}]P(x_1,\ldots,x_n)\not=0$. Let $A_1,\ldots,A_n$ be finite
subsets of $F$ with $|A_i|>k_i$ for $i=1,\ldots,n$. Then, for the
restricted sumset
\begin{equation}\label{1.4}C=\{x_1+\cdots+x_n:\ x_1\in A_1,\ldots,x_n\in A_n,\ \t{and}\ P(x_1,\ldots,x_n)\not=0\},
\end{equation}
we have
\begin{equation}\label{1.5}|C|\gs\min\{p(F)-\deg P,\ |A_1|+\cdots+|A_n|-n-2\deg P+1\}.\end{equation}
\end{theorem}

\begin{remark} Theorem \ref{Th1.3} in the case $P(x_1,\ldots,x_n)=1$
gives the Cauchy-Davenport theorem. When $F$ is of characteristic
zero (i.e., $p(F)=+\infty$), Theorem \ref{Th1.3} extends a result of Sun
\cite[Theorem 1.1]{S01} on sums of subsets of $\Z$ with various linear
restrictions.
\end{remark}

The following example shows that the lower bound in Theorem \ref{Th1.3} is
essentially best possible.

\medskip

\noindent{\it Example} 1.2. Let $p$ be a prime and let $F_p$ be the
finite field $\Z/p\Z$.

{(i) Let $$P(x_1,\ldots,x_n)=\prod_{s\in S}(x_1+\cdots+x_n-s)$$
where $S$ is a nonempty subset of $F_p$. Then
\begin{align*}&|\{x_1+\cdots+x_n:\ x_1,\ldots,x_n\in F_p\ \t{and}\
P(x_1,\ldots,x_n)\not=0\}|
\\&\qquad =|F_p\sm S|=|F_p|-|S|=p-\deg P.
\end{align*}

{(ii)} Let $A=\{\bar r=r+p\Z:\ r\in[0,m-1]\}\se F_p$ with $n\ls m\ls
p$, where $n$ is a positive integer. If $p\gs n(m-n)+1$, then
\begin{align*}&|\{x_1+\cdots+x_n:\,
x_1,\ldots,x_n\in A,\ \t{and}\ x_i\not=x_j\ \t{if}\ i\not=j\}|
\\=&|\{\bar r:\ r\in [0+\cdots+(n-1),\, (m-n)+\cdots+(m-1)]\}|
\\=&n(m-n)+1=n|A|-n-2\deg\prod_{1\ls i<j\ls n}(x_j-x_i)+1.
\end{align*}

\medskip

Here are some consequences of Theorem \ref{Th1.3}.

\begin{corollary}\label{Cor1.2} Let $A$ be a finite subset of a field $F$,
and let $a_1,\ldots,a_n\in F^*$.

{\rm (i)} For any $f(x)\in F[x]$ with $\deg f=m\gs0$, we have
\begin{equation}\label{1.6}\begin{aligned}&|\{a_1x_1+\cdots+a_nx_n:\ x_1,\ldots,x_n\in A,\ \t{and}\ f(x_i)\not=f(x_j)\ \t{if}\ i\not=j\}|
\\&\qquad\gs\min\l\{p(F)-m\bi n2,\ n(|A|-1-m(n-1))+1\r\}.
\end{aligned}\end{equation}

{\rm (ii)} Let $S_{ij}\se F$ with $|S_{ij}|\ls2m-1$ for all $1\ls
i<j\ls n$. Then
\begin{equation}\label{1.7}\begin{aligned}&|\{a_1x_1+\cdots+a_nx_n:\ x_1,\ldots,x_n\in A,\ \t{and}\ x_i-x_j\not\in S_{ij}\ \t{if}\ i<j\}|
\\&\gs\min\l\{p(F)-(2m-1)\bi n2,\ n(|A|-1-(2m-1)(n-1))+1\r\}.
\end{aligned}\end{equation}
\end{corollary}

\begin{remark} In the case $m=1$, each of the two parts in Corollary
1.2 yields the inequality
\begin{equation}\label{1.8}\begin{aligned}&|\{a_1x_1+\cdots+a_nx_n:\ x_1,\ldots,x_n\in A,\ \t{and}\ x_i\not=x_j\ \t{if}\ i\not=j\}|
\\&\qquad \gs\min\l\{p(F)-\bi n2,\ n(|A|-n)+1\r\}.
\end{aligned}\end{equation}
\end{remark}

Let $m_1,\cs,m_n\in\N$. When we expand $\prod_{1\le i,j\le n,\
i\not=j}(1-x_i/x_j)^{m_j}$ as a Laurent polynomial (with negative
exponents allowed), the constant term was conjectured to be the
multinomial coefficient $(\sum_{i=1}^nm_i)!/\prod_{i=1}^nm_i!$ by F.
J. Dyson \cite{D62} in 1962. A simple proof of Dyson's conjecture given
by I. J. Good \cite{G70} employs the Lagrange interpolation formula.
Using Dyson's conjecture we can deduce the following result from
Theorem 1.3.

\begin{corollary}\label{Cor1.3} Let $A_1,\ldots,A_n\ (n>1)$ be finite
nonempty subsets of a field $F$, and let $S_{ij}\ (1\ls i\not=j\ls
n)$ be subsets of $F$ with $|S_{ij}|\ls (|A_i|-1)/(n-1)$. Then, for
any $a_1,\ldots,a_n\in F^*$, we have
\begin{equation}\begin{aligned}\label{1.9}&|\{a_1x_1+\cdots+a_nx_n:\ x_1\in A_1,\ldots,x_n\in A_n,\ \t{and}\ x_i-x_j\not\in S_{ij}\ \t{if}\ i\not=j\}|
\\&\gs\min\l\{p(F)-(n-1)\sum_{i=1}^nm_i,\ \sum_{i=1}^n(|A_i|-1)-2(n-1)\sum_{i=1}^n m_i+1\r\},
\end{aligned}\end{equation}
where $m_i=\max_{j\in[1,n]\sm\{i\}}|S_{ij}|$ for $i=1,\ldots,n$.
\end{corollary}

In the next section we will prove Theorem \ref{Th1.2} with the help of
several lemmas. Section 3 is devoted to the proof of Theorem \ref{Th1.3}.
Theorem \ref{Th1.1} and Corollaries 1.1-1.3 will be shown in Section 4.
Finally, in Section 5 we deduce a further extension of Theorem \ref{Th1.3}.

\section{Proof of Theorem \ref{Th1.2}}

\setcounter{lemma}{0}
\setcounter{theorem}{0}
\setcounter{corollary}{0}
\setcounter{remark}{0}
\setcounter{equation}{0}
\setcounter{conjecture}{0}

\begin{lemma}\label{Lem2.1} Let $a_1,\ldots,a_n$ be nonzero elements in a
field $F$ with $p(F)\not=2$. Then, for some $\sigma\in S_n$ we have
$$a_{\sigma(2i-1)}+a_{\sigma(2i)}\not=0\ \quad\t{for all}\
0<i\ls\l\lfloor\f n2\r\rfloor-\da(a_1,\ldots,a_n),$$
where $\da(a_1,\ldots,a_n)\in\{0,1\}$ takes the value $1$ if and only if there exists $a\in F^*$ such that
$\{a_1,\ldots,a_n\}=\{a,-a\}$ and
\begin{equation}\label{2.1}|\{1\ls i\ls n:\ a_i=a\}|\eq|\{1\ls i\ls n:\ a_i=-a\}|\eq1\ (\mo\ 2).\end{equation}
\end{lemma}
\Proof. We use induction on $n$.

The case $n\in\{1,2\}$ is trivial.

Now let $n>2$ and assume the desired result for smaller values of $n$.
\medskip

In the case $\da(a_1,\ldots,a_n)=1$, there is an element $a\in F^*$ such that $\{a_1,\ldots,a_n\}=\{a,-a\}$
and (2.1) holds; thus the desired result follows immediately since $a+a\not=0$ and $-a+(-a)\not=0$.
\medskip

Below we let $\da(a_1,\ldots,a_n)=0$. If $a_1+a_2=a_1+a_3=a_2+a_3=0$, then $a_1=a_2=a_3=0$
which contradicts the condition $a_1,\ldots,a_n\in F^*$.
So for some $1\ls s<t\ls n$ we have $a_s+a_t\not=0$. Without loss of generality
we simply suppose that $a_{n-1}+a_n\not=0$.
By the induction hypothesis, for some $\sigma\in S_{n-2}$ we have
$$a_{\sigma(2i-1)}+a_{\sigma(2i)}\not=0\ \quad\t{for all}\
\ 0<i\ls\l\lfloor\f {n-2}2\r\rfloor-\da(a_1,\ldots,a_{n-2}).$$

If $\da(a_1,\ldots,a_{n-2})=0$, then it suffices to set $\sigma(2\lfloor n/2\rfloor-1)=n-1$
and $\sigma(2\lfloor n/2\rfloor)=n$.

Now let $\da(a_1,\ldots,a_{n-2})=1$. Then for some $a\in F^*$ we have both
$\{a_1,\ldots,a_{n-2}\}=\{a,-a\}$ and
$$|\{1\ls i\ls n-2:\ a_i=a\}|\eq|\{1\ls i\ls n-2:\ a_i=-a\}|\eq1\ (\mo\ 2).$$

\medskip

{\it Case} 1. $\{a,-a\}\cap\{a_{n-1},a_n\}=\emptyset$.

 In this case, $a+a_{n-1}\not=0$ and $-a+a_n\not=0$.
Thus there exists $\sigma\in S_n$ such that $a_{\sigma(2i-1)}=a_{\sigma(2i)}\in\{a,-a\}$ for
all $0<i<\lfloor (n-2)/2\rfloor$, and also
$$a_{\sigma(2\lfloor (n-2)/2\rfloor-1)}=a,\ \ a_{\sigma(2\lfloor (n-2)/2\rfloor)}=a_{n-1}$$
and
$$a_{\sigma(2\lfloor n/2\rfloor-1)}=-a,\ \ a_{\sigma(2\lfloor n/2\rfloor)}=a_{n}.$$

\medskip

{\it Case} 2. $\{a,-a\}\cap\{a_{n-1},a_n\}\not=\emptyset$.

Without loss of generality we assume that $a_{n-1}=a$. As $\da(a_1,\ldots,a_n)=0$ we cannot have $a_{n-1}=a_n\in\{a,-a\}$.
Thus $a_n\not=a$. Now $a+a_{n-1}=2a\not=0$ and $-a+a_n\not=0$. As in Case 1 there exists $\sigma\in S_n$ such that
$a_{\sigma(2i-1)}=a_{\sigma(2i)}\in\{a,-a\}$ for
all $0<i\ls\lfloor n/2\rfloor$.

So far we have proved the desired result by induction. \qed

\begin{lemma}\label{Lem2.2} Let $k_1,\ldots,k_n\in\N$ and
$a_1,\ldots,a_n\in F^*$, where $F$ is a field with $p(F)\not=2$. Set
\begin{equation}\label{2.2}
f(x_1,\ldots,x_n)=\sum_{\sigma\in S_n}\sgn(\sigma)\prod_{j=1}^n(k_j-x_j)_{\sigma(j)-1}\,a_j^{\sigma(j)-1}
\end{equation}
and let $\da(a_1,\ldots,a_n)$ be as in Lemma 2.1.
Provided the following {\rm (i)} or {\rm (ii)},
there are $m_1,\ldots,m_n\in\N$ not exceeding $\max\{2n-3,0\}$ such that $m_1+\cdots+m_n=\bi n2$
and $f(m_1,\ldots,m_n)\not=0$.

{\rm (i)} $\da(a_1,\ldots,a_n)=0$.

{\rm (ii)} $\da(a_1,\ldots,a_n)=1$, and for some $1\ls s<t\ls n$ we have $a_s+a_t=0$ and $k_s+k_t\not\eq1\ (\mo\ p(F))$.
(A congruence modulo $\infty$ refers to the corresponding equality.)
\end{lemma}
\Proof. We use induction on $n$.

 When $n=1$, obviously we can take $m_1=\cdots=m_n=0$ to meet the requirement.

 In the case $n=2$, we have $f(x_1,x_2)=a_2(k_2-x_2)-a_1(k_1-x_1)$. Clearly $f(1,0)-f(0,1)=a_1+a_2$. If $f(1,0)=f(0,1)$, then
 $a_1+a_2=0$, $\da(a_1,a_2)=1$ and $f(0,1)=a_2(k_2-1)-a_1k_1=a_2(k_1+k_2-1)\not=0$ by condition (ii).
 Anyway, we have $f(m_1,m_2)\not=0$ for some $m_1\in\{0,1\}$ and $m_2=1-m_1$.

 Below we let $n\gs3$ and assume the desired result for smaller values of $n$.
In case (ii), clearly $\da(a_3,\ldots,a_n)=0$, and we may simply assume that $s=1$ and $t=2$ without loss of generality. By Lemma \ref{Lem2.1}, there
is a rearrangement $a_1',\ldots,a_n'$  of $a_1,\ldots,a_n$
 such that $a_{n-2i-1}'+a_{n-2i}'\not=0$ for all
$0\ls i<\lfloor n/2\rfloor-\da(a_1,\ldots,a_n)$, and $a_1'=a_1$ and $a_2'=a_2$ in case (ii).
Suppose that $a_i'=a_{\tau(i)}$ for $i=1,\ldots,n$, where $\tau\in S_n$, and $\tau(1)=1$ and $\tau(2)=2$ in case (ii).
Set $k_i'=k_{\tau(i)}$ for $i=1,\ldots,n$. Then
\begin{align*} f(x_1,\ldots,x_n)=&\sum_{\sigma\in S_n}\sgn(\sigma)
\prod_{i=1}^n(k_{\tau(i)}-x_{\tau(i)})_{\sigma\tau(i)-1}(a_i')^{\sigma\tau(i)-1}
\\=&\sgn(\tau)\sum_{\pi\in S_n}\sgn(\pi)\prod_{i=1}^n(k_i'-x_{\tau(i)})_{\pi(i)-1}(a_i')^{\pi(i)-1}
\end{align*}
Hence $f(m_1,\ldots,m_n)\not=0$ for some $m_1,\ldots,m_n\in[0,2n-3]$ if and only if
$$\sum_{\pi\in S_n}\sgn(\pi)\prod_{i=1}^n(k_i'-m_i')_{\pi(i)-1}(a_i')^{\pi(i)-1}\not=0$$
for some $m_1',\ldots,m_n'\in[0,2n-3]$. Without loss of generality, below we simply assume that $a_i'=a_i$
and $k_i'=k_i$ for all $i=1,\ldots,n$.

By the induction hypothesis, there are $m_1,\ldots,m_{n-2}\in[0,2n-3]$ such that $\sum_{j=1}^{n-2} m_j=\bi{n-2}2$
and
$$\Sigma:=\sum_{\sigma\in S_{n-2}}\sgn(\sigma)\prod_{j=1}^{n-2}(k_j-m_j)_{\sigma(j)-1}\,a_j^{\sigma(j)-1}\not=0$$
Define
\begin{align*} g(x)=&f\l(m_1,\ldots,m_{n-2},x,\bi n2-x-m_1-\cdots-m_{n-2}\r)
\\=&\sum_{\sigma\in S_n}\sgn(\sigma)\prod_{j=1}^{n-2}(k_j-m_j)_{\sigma(j)-1}\,a_j^{\sigma(j)-1}
\\&\ \ \times (k_{n-1}-x)_{\sigma(n-1)-1}\,a_{n-1}^{\sigma(n-1)-1}
\\&\ \ \times \bigg(k_n-\bi n2+x+\sum_{j=1}^{n-2}m_j\bigg)_{\sigma(n)-1}\,a_n^{\sigma(n)-1}.
\end{align*}
For $\sigma\in S_n$, if $\sigma(1)-1+(\sigma(2)-1)=2n-3$ then $\{\sigma(1),\sigma(2)\}=\{n-1,n\}$.
Thus
\begin{align*}[x^{2n-3}]g(x)=&\sum_{\sigma\in S_n\atop\{\sigma(n-1),\sigma(n)\}=\{n-1,n\}}\sgn(\sigma)
\prod_{j=1}^{n-2}(k_j-m_j)_{\sigma(j)-1}\,a_j^{\sigma(j)-1}
\\&\qquad\qquad\qquad\times(-a_{n-1})^{\sigma(n-1)-1}a_n^{\sigma(n)-1}
\\=&\sum_{\sigma\in S_{n-2}}\sgn(\sigma)
\prod_{j=1}^{n-2}(k_j-m_j)_{\sigma(j)-1}\,a_j^{\sigma(j)-1}
\\&\quad \times\l((-a_{n-1})^{n-2}a_n^{n-1}-(-a_{n-1})^{n-1}a_n^{n-2}\r)
\\=&(-1)^n(a_{n-1}a_n)^{n-2}(a_{n-1}+a_n)\Sigma\not=0.
\end{align*}
Since $\deg g(x)=2n-3$, there is an integer $m_{n-1}\in[0,2n-3]$ such that $g(m_{n-1})\not=0$.
Set
$$m_n=\bi n2-\sum_{j=1}^{n-1}m_j=\bi n2-\bi{n-2}2-m_{n-1}=2n-3-m_{n-1}.$$ Then
$$f(m_1,\ldots,m_n)=g(m_{n-1})\not=0.$$
This concludes the induction step and we are done. \qed

\begin{lemma}\label{Lem2.3} Let $F$ be a field with $p(F)\not=2$, and let $a_1,\ldots,a_n\ (n\gs4)$
be nonzero elements of $F$ with $\da(a_1,\ldots,a_n)=1$.
Suppose that $p(F)\gs\sum_{j=1}^nk_j-n^2+n+1$ where $k_1,\ldots,k_n$ are integers not smaller than $2n-3$.
Then there are $1\ls s<t\ls n$ such that $a_s+a_t=0$ and $k_s+k_t\not\eq1\ (\mo\ p(F))$,
unless $n=4$ and there is a permutation
$\sigma\in S_4$ such that $a_{\sigma(1)}=a_{\sigma(2)}=a_{\sigma(3)}$, $k_{\sigma(1)}=k_{\sigma(2)}=k_{\sigma(3)}=5$
and $k_{\sigma(4)}=p(F)-4$.
\end{lemma}
\Proof. For any $1\ls s<t\ls n$ we have
\begin{align*} p(F)-(k_s+k_t-1)\gs&\sum_{1\ls j\ls n\atop j\not=s,t} k_j-n^2+n+2
\\\gs&(n-2)(2n-3)-n^2+n+2=(n-2)(n-4)
\end{align*}
and hence
\begin{align*} &k_s+k_t\eq 1\ (\mo\ p(F))
\\\iff& k_s+k_t-1=p(F),\ k_i=2n-3\ \t{for}\ i\in[1,n]\sm\{s,t\},\ \t{and}\ n=4.
\end{align*}
Since $\da(a_1,\ldots,a_n)=1$, for some $1\ls s<t\ls n$ we have
$a_s+a_t=0$; also $k_s+k_t\not\eq1\ (\mo\ p(F))$ if $n>4$. This
proves the desired result for $n>4$.

Now assume $n=4$. By $\da(a_1,a_2,a_3,a_4)=1$, there is a permutation $\sigma\in S_4$ such that
$a_{\sigma(1)}=a_{\sigma(2)}=a_{\sigma(3)}=-a_{\sigma(4)}$. Clearly $a_{\sigma(i)}+a_{\sigma(4)}=0$ for any $i=1,2,3$.
Suppose that $k_{\sigma(i)}+k_{\sigma(4)}\eq1\ (\mo\ p(F))$ for all $i=1,2,3$. By the above,
$k_{\sigma(i)}+k_{\sigma(4)}-1=p(F)$ for $i=1,2,3$, and $k_{\sigma(1)}=k_{\sigma(2)}=k_{\sigma(3)}=2n-3=5$.
It follows that $k_{\sigma(4)}=p(F)-4$.

The proof of Lemma \ref{Lem2.3} is now complete. \qed

 \begin{lemma}\label{Lem2.4} Let $F$ be a field of prime characteristic with $p(F)=p\gs7$ and let $a_1=a_2=a_3=a\in F^*$ and $a_4=-a$.
 Let $k_1=k_2=k_3=5$ and $k_4=p-4$.
 Then there are $m_1,m_2,m_3,m_4\in[0,3]$ such that $m_1+m_2+m_3+m_4=\bi 42=6$ and
 $$\sum_{\sigma\in S_4}\sgn(\sigma)\prod_{j=1}^4(k_j-m_j)_{\sigma(j)-1}\,a_j^{\sigma(j)-1}\not=0.$$
 \end{lemma}
 \Proof. Set $m_1=0$, $m_2=2$, $m_3=3$ and $m_4=1$. Then
 \begin{align*} &\sum_{\sigma\in S_4}\sgn(\sigma)\prod_{j=1}^4(k_j-m_j)_{\sigma(j)-1}\,a_j^{\sigma(j)-1}
\\ =&\sum_{\sigma\in S_4}\sgn(\sigma)\prod_{j=1}^3(5-m_j)_{\sigma(j)-1}
\times(-4-m_4)_{\sigma(4)-1}(-1)^{\sigma(4)-1}a^{0+1+2+3}
\\=&-480a^6\not=0
\end{align*}
since $p$ does not divide $480$. We are done. \qed

\medskip

\noindent {\it Proof of Theorem \ref{Th1.2}}.
 Set $A_i'=a_iA_i=\{a_ix_i:\ x_i\in A_i\}$ and $a_i'=a_i^{-1}$ for $i=1,\ldots,n$. Then
$$C=\{y_1+\cdots+y_n:\ y_1\in A_1',\ldots,y_n\in A_n',\ \t{and}\ a_i'y_i\not=a_j'y_j\ \t{if}\ i\not=j\}.$$
In the case $n=1$, clearly
$$|C|=|A_1'|=|A_1|\gs\min\{p(F),\,|A_1|-1^2+1\}.$$
When $n=2$, we have
\begin{align*}|C|=&|\{y_1+y_2:\, y_1\in A_1',\ y_2\in A_2'\ \t{and}\ y_1-(a_1')^{-1}a_2'y_2\not=0\}|
\\\gs&\min\{p(F)-[\![a_1'=a_2']\!],\,|A_1'|+|A_2'|-3\}\ \ \t{(by \cite[Corollary 3]{PS1})}
\\=&\min\{p(F)-[\![a_1=a_2]\!],\,|A_1|+|A_2|-2^2+1\}.
\end{align*}

Below we let $n>2$. Clearly $p(F)\gs(n-1)^2>2$. Define
\begin{equation}\label{2.3}N=\sum_{j=1}^n|A_j|-n^2.\end{equation}
We want to show that $|C|\gs\min\{p(F),\,N+1\}$.

Let's first assume that $p(F)>N$. Note that $p(F)\gs(4-1)^2>7$ if
$n\gs 4$. In view of Lemmas 2.1-2.4, there are
$m_1,\ldots,m_n\in[0,2n-3]$ such that $m_1+\cdots+m_n=\bi n2$ and
\begin{equation}\label{2.4}
S=\sum_{\sigma\in S_n}\sgn(\sigma)\prod_{j=1}^n(|A_j'|-1-m_j)_{\sigma(j)-1}(a_j')^{\sigma(j)-1}\not=0.
\end{equation}
Clearly it suffices to deduce a contradiction under the assumption that $|C|\ls N$.
Let $P(x_1,\ldots,x_n)$ be the polynomial
$$\prod_{1\ls i<j\ls n}(a_j'x_j-a_i'x_i)\times\prod_{j=1}^nx_j^{m_j}
\times \prod_{x\in C}(x_1+\cdots+x_n-c)\times
(x_1+\cdots+x_n)^{N-|C|}$$ Then $\deg P=\sum_{j=1}^n(|A_j'|-1)$,
since
\begin{align*}&[x_1^{|A_1'|-1}\cdots x_n^{|A_n'|-1}]P(x_1,\ldots,x_n)
\\=&\bigg[\prod_{j=1}^nx_j^{|A_j'|-1-m_j}\bigg]\sum_{\sigma\in S_n}\sgn(\sigma)\prod_{j=1}^n(a_j'x_j)^{\sigma(j)-1}
\times(x_1+\cdots+x_n)^N
\\=&\sum_{\sigma\in S_n\atop\sigma(j)\ls|A_j'|-m_j\ \t{for}\ j\in[1,n]}
\sgn(\sigma)\f{N!}{\prod_{j=1}^n(|A_j'|-m_j-\sigma(j))!}\prod_{j=1}^n(a_j')^{\sigma(j)-1}
\end{align*}
and hence
$$\prod_{j=1}^n(|A_j'|-1-m_j)!\times [x_1^{|A_1'|-1}\cdots x_n^{|A_n'|-1}]P(x_1,\ldots,x_n)=N!S\not=0.$$
Thus, by the Combinatorial Nullstellensatz there are $y_1\in A_1',\ldots,y_n\in A_n'$ such that
$P(y_1,\ldots,y_n)\not=0$ which contradicts the definition of $C$.

\medskip

Now we handle the case $p(F)\ls N$. Since $n(2n-2)-n^2\ls p(F)-1<\sum_{j=1}^n|A_j|-n^2$,
we can choose $B_j\se A_j$ with $|B_j|\gs 2n-2$ so that $M=\sum_{j=1}^n|B_j|-n^2=p(F)-1$.
As $p(F)>M$, by the above we have
\begin{align*}|C|\gs &|\{a_1x_1+\cdots+a_nx_n:\ x_1\in B_1,\ldots,x_n\in B_n,\ \t{and}\ x_i\not=x_j\ \t{if}\ i\not=j\}|
\\\gs& M+1=\min\{p(F),N\}.
\end{align*}

The proof of Theorem \ref{Th1.2} is now complete. \qed

\section{Proof of Theorem \ref{Th1.3}}

\setcounter{lemma}{0}
\setcounter{theorem}{0}
\setcounter{corollary}{0}
\setcounter{remark}{0}
\setcounter{equation}{0}
\setcounter{conjecture}{0}

The inequality (\ref{1.5}) holds trivially if $p(F)\ls \deg P$ or $\sum_{i=1}^n|A_i|<n+2\deg P$.
Below we assume that $p(F)>\deg P$ and $\sum_{i=1}^n|A_i|\gs n+2\deg P.$

Write
\begin{equation}\label{3.1}P(x_1,\ldots,x_n)=\sum_{j_1,\ldots,j_n\gs0\atop\j_1+\cdots+j_n\ls\deg P}
c_{j_1,\ldots,j_n}x_1^{j_1}\cdots x_n^{j_n}
\ \ \t{with}\ c_{j_1,\ldots,j_n}\in F,
\end{equation}
and define
\begin{equation}\label{3.2}
P^*(x_1,\ldots,x_n)=\sum_{j_1,\ldots,j_n\gs0\atop\j_1+\cdots+j_n=\deg P}
c_{j_1,\ldots,j_n}(x_1)_{j_1}\cdots(x_n)_{j_n} \in
F[x_1,\ldots,x_n].
\end{equation}
It is easy to see that
 $$[x_1^{k_1}\cdots x_n^{k_n}]P^*(x_1,\ldots,x_n)
=[x_1^{k_1}\cdots x_n^{k_n}]P(x_1,\ldots,x_n)\not=0.$$

 To distinguish from the integer 1, we use $e$ to denote the multiplicative identity of the field $F$.
For each $i=1,\ldots,n$, clearly the set
$$B_i=\{me:\ m\in[|A_i|-k_i-1,|A_i|-1]\}$$
has cardinality $k_i+1$ since $k_i\ls \deg P<p(F)$.
Thus, by the Combinatorial Nullstellensatz, there are
\begin{equation}\label{3.3}m_1\in[|A_1|-k_1-1,|A_1|-1],\ \ldots,\ m_n\in[|A_n|-k_n-1,|A_n|-1]
\end{equation}
such that
\begin{equation}\label{3.4}P^*(m_1e,\ldots,m_ne)\not=0.\end{equation}

Define
\begin{equation}\label{3.5}M=m_1+\cdots+m_n-\deg P.\end{equation}
Clearly $$M\gs\sum_{i=1}^n(|A_i|-k_i-1)-\deg P=\sum_{i=1}^n|A_i|-n-2\deg P\gs0.$$
Observe that
  \begin{align*} &[x_1^{m_1}\cdots x_n^{m_n}]P(x_1,\ldots,x_n)(x_1+\cdots+x_n)^M
 \\=&\sum_{j_1\in[0,m_1],\ldots,j_n\in[0,m_n]\atop\j_1+\cdots+j_n=\deg P}
 \f{M!}{(m_1-j_1)!\cdots(m_n-j_n)!}c_{j_1,\ldots,j_n}
 \end{align*}
 and thus
\begin{align*}&m_1!\cdots m_n![x_1^{m_1}\cdots x_n^{m_n}]P(x_1,\ldots,x_n)(x_1+\cdots+x_n)^M
 \\=&M!\sum_{j_1,\ldots,j_n\gs0\atop\j_1+\cdots+j_n=\deg P} (m_1e)_{j_1}\cdots(m_ne)_{j_n}c_{j_1,\ldots,j_n}
 =M!P^*(m_1e,\ldots,m_ne).
 \end{align*}
In the case $|C|\ls M<p(F)$,  with the help of (3.4) we have
\begin{align*} &[x_1^{m_1}\cdots x_n^{m_n}]P(x_1,\ldots,x_n)(x_1+\cdots+x_n)^{M-|C|}\prod_{c\in C}(x_1+\cdots+x_n-c)
\\&\quad=[x_1^{m_1}\cdots x_n^{m_n}]P(x_1,\ldots,x_n)(x_1+\cdots+x_n)^M\not=0,
\end{align*}
hence by the Combinatorial Nullstellensatz there are $x_1\in A_1,\ldots,x_n\in A_n$ such that
$$P(x_1,\ldots,x_n)(x_1+\cdots+x_n)^{M-|C|}\prod_{c\in C}(x_1+\cdots+x_n-c)\not=0$$
which is impossible by the definition of $C$. Therefore, either
\begin{equation}\label{3.6}p(F)\ls M\ls \sum_{i=1}^n(|A_i|-1)-\deg P\end{equation}
 or
 \begin{equation}\label{3.7}|C|\gs M+1\gs \sum_{i=1}^n|A_i|-n-2\deg P+1.\end{equation}

If $p(F)>\sum_{i=1}^n(|A_i|-1)-\deg P$, then (\ref{3.6}) fails and hence
\begin{align*}|C|\gs& \sum_{i=1}^n|A_i|-n-2\deg P+1
\\=&\min\bg\{p(F)-\deg P,\ \sum_{i=1}^n|A_i|-n-2\deg P+1\bg\}.
\end{align*}

In the case $p(F)\ls \sum_{i=1}^n(|A_i|-1)-\deg P$, as $\sum_{i=1}^n k_i=\deg P$
there are $A_1'\se A_1,\ldots,A_n'\se A_n$ such that
$$|A_1'|>k_1,\ldots,|A_n'|>k_n,\ \t{and}\ \sum_{i=1}^n(|A_i'|-1)-\deg P=p(F)-1<p(F),$$
therefore
\begin{align*} |C|\gs&|\{x_1+\cdots+x_n:\ x_1\in A_1',\ldots,x_n\in A_n',\ \t{and}\ P(x_1,\ldots,x_n)\not=0\}|
\\\gs&\min\bg\{p(F)-\deg P,\ \sum_{i=1}^n|A_i'|-n-2\deg P+1\bg\}
\\=&p(F)-\deg P=\min\bg\{p(F)-\deg P,\ \sum_{i=1}^n|A_i|-n-2\deg P+1\bg\}.
\end{align*}
This concludes the proof.

\section{Proofs of Corollaries 1.1-1.3 and Theorem \ref{Th1.1}}

\setcounter{lemma}{0}
\setcounter{theorem}{0}
\setcounter{corollary}{0}
\setcounter{remark}{0}
\setcounter{equation}{0}
\setcounter{conjecture}{0}

\noindent{\it Proof of Corollary \ref{Cor1.1}}. As $A$ has a subset of cardinality $\lceil \sqrt{4p-7}\ \rceil$,
it suffices to consider the case $|A|=\lceil \sqrt{4p-7}\ \rceil$.
Since $n-1\ls |A|/2-1<\sqrt p$ and $(n-A|/2)^2\ls |A|^2/4-p+1$, applying Theorem \ref{Th1.2} we get
$$\{a_1x_1+\cdots+a_nx_n:\ x_1,\ldots,x_n\in A\ \t{and}\ x_i\not=x_j\ \t{if}\ i\not=j\}=F_p.$$
This concludes the proof. \qed

\medskip

\noindent{\it Proof of Corollary \ref{Cor1.2}}. Both (\ref{1.6}) and (\ref{1.7}) are trivial in the case $|A|\ls m(n-1)$.
Below we assume that $|A|>m(n-1)$, and put $A_i=\{a_ix:\, x\in A\}$ for $i=1,\ldots,n$.

(i) Set $b_j=[x^m]f(x)a_j^{-m}$ for $j\in[1,n]$, and define
$$P(x_1,\ldots,x_n)=\prod_{1\ls i<j\ls n}\l(f(a_j^{-1}x_j)-f(a_i^{-1}x_i)\r).$$
Note that
\begin{align*}\prod_{1\ls i<j\ls n}(b_jx_j^m-b_ix_i^m)
=&\det((b_jx_j^m)^{i-1})_{1\ls i,j\ls n}
\\=&\sum_{\sigma\in S_n}\sgn(\sigma)
\prod_{i=1}^nb_{\sigma(i)}^{i-1}x_{\sigma(i)}^{(i-1)m}.
\end{align*}
Therefore
$$\bigg[\prod_{i=1}^nx_i^{(i-1)m}\bigg]P(x_1,\ldots,x_n)\not=0\ \ \t{and}\ \ \sum_{i=1}^n(i-1)m=\deg P.$$

 In view of Theorem \ref{Th1.3},
$$\begin{aligned}&|\{a_1x_1+\cdots+a_nx_n:\ x_1,\ldots,x_n\in A,\ \t{and}\ f(x_i)\not=f(x_j)\ \t{if}\ i\not=j\}|
\\=&|\{y_1+\cdots+y_n:\ y_1\in A_1,\ldots,y_n\in A_n,\ \t{and}\ P(y_1,\ldots,y_n)\not=0\}|
\\\gs&\min\l\{p(F)-\deg P,\ |A_1|+\cdots+|A_n|-n-2\deg P+1\r\}
\\=&\min\l\{p(F)-m\bi n2,\ n(|A|-1)-mn(n-1)+1\r\}.
\end{aligned}$$
So we have (\ref{1.6}).

(ii) Let $P(x_1,\ldots,x_n)$ be the polynomial
$$\prod_{1\ls i<j\ls n}\((a_j^{-1}x_j-a_i^{-1}x_i)^{2m-1-|S_{ij}|}
\prod_{s\in S_{ij}}(a_j^{-1}x_j-a_i^{-1}x_i+s)\).$$
By \cite[(2.8)]{SY},
\begin{align*}&\bigg[\prod_{i=1}^nx_i^{(m-1)(n-1)+i-1}\bigg]P(a_1x_1,\ldots,a_nx_n)
\\=&\bigg[\prod_{i=1}^nx_i^{(m-1)(n-1)+i-1}\bigg]\prod_{1\ls i<j\ls n}(x_j-x_i)^{2m-1}
=(-1)^{(m-1)\bi n2}Ne,
\end{align*}
where $N=(mn)!/(m!^nn!)\in\Z^+=\{1,2,3,\ldots\}$. Clearly $N=1$ if $m=1$ or $n=1$.
 If $\min\{m,n\}\gs2$ and $mn\gs p(F)$, then
\begin{align*} p(F)-(2m-1)\bi n2\ls &mn-1-\l(m-\f12\r)n(n-1)
\\=&n\l(m-\l(m-\f12\r)(n-1)\r)-1\ls0.\end{align*}
So (\ref{1.7}) holds trivially if $mn\gs p(F)$.

 Below we handle the case $mn<p(F)$, thus $Ne\not=0$. Note that
$$\bigg[\prod_{i=1}^nx_i^{(m-1)(n-1)+i-1}\bigg]P(x_1,\ldots,x_n)\not=0.$$
Clearly $\sum_{i=1}^n((m-1)(n-1)+i-1)=(2m-1)\bi n2=\deg P$.
Observe that $|A_i|=|A|>m(n-1)\gs(m-1)(n-1)+i-1$ for all $i\in[1,n]$.
Applying Theorem \ref{Th1.3} we get
$$\begin{aligned}&|\{a_1x_1+\cdots+a_nx_n:\ x_1,\ldots,x_n\in A,\ \t{and}\ x_i-x_j\not\in S_{ij}\ \t{if}\ i<j\}|
\\\gs&|\{y_1+\cdots+y_n:\ y_1\in A_1,\ldots,y_n\in A_n,\ \t{and}\ P(y_1,\ldots,y_n)\not=0\}|
\\\gs&\min\l\{p(F)-\deg P,\ |A_1|+\cdots+|A_n|-n-2\deg P+1\r\}
\\=&\min\l\{p(F)-(2m-1)\bi n2,\ n(|A|-1)-(2m-1)n(n-1)+1\r\}.
\end{aligned}$$
This proves (\ref{1.7}).

So far we have completed the proof of Corollary \ref{Cor1.2}. \qed

\medskip

The Dyson conjecture mentioned in Section 1 can be restated as follows: For any $m_1,\ldots,m_n\in\N$ we have
\begin{equation}\label{4.1}\begin{aligned}&[x_1^{m_1(n-1)} \cdots x_n^{m_n(n-1)}]\prod_{1 \le i < j \le n} (x_i - x_j)^{m_i + m_j}
\\&\qquad=(-1)^{\sum_{j=1}^n (j-1)m_j} \f{(m_1 + \cdots + m_n)!}{m_1! \cdots m_n!}.
\end{aligned}\end{equation}
A combinatorial proof of this was given by D. Zeilberger \cite{Z82} in 1982.
Below we use (\ref{4.1}) to prove Corollary \ref{Cor1.3}.
\medskip

\noindent{\it Proof of Corollary \ref{Cor1.3}}. We only need to consider the nontrivial case
 $\sum_{i=1}^n m_i<p(F)$.  Similar to the proof of Corollary \ref{Cor1.2},
it suffices to note that the coefficient of the monomial
$\prod_{i=1}^nx_i^{m_i(n-1)}$ in the polynomial $\prod_{1 \le i < j
\le n} (x_i - x_j)^{m_i + m_j}$ over $F$ does not vanish by (\ref{4.1})
and $\sum_{k=1}^nm_k<p(F)$. \qed

\medskip
\noindent{\it Proof of Theorem \ref{Th1.1}}. If $p(F)-\bi n2\gs n|A|-n^2+1$, then (\ref{1.1}) follows from (\ref{1.8}).

Now assume that $p(F)-\bi n2\ls n|A|-n^2$. Then
$$n|A|\gs p(F)-\bi n2+n^2\gs\f{3n^2-5n}2-\f{n^2-n}2+n^2=2n^2-2n$$
and hence $|A|\gs 2n-2$. Note also that if $n>1$ then $p(F)\gs n(3n-5)/2\gs (n-1)^2$.
Thus, by applying Theorem \ref{Th1.2} we obtain the desired result. \qed

\section{A Further Extension of Theorem \ref{Th1.3}}

\setcounter{lemma}{0}
\setcounter{theorem}{0}
\setcounter{corollary}{0}
\setcounter{remark}{0}
\setcounter{equation}{0}
\setcounter{conjecture}{0}

Recently Z. W. Sun \cite{S08a} employed the Combinatorial Nullstellensatz
to establish the following result on value sets of polynomials.

\begin{theorem}[Sun \cite{S08a}]\label{Th5.1}  Let
$A_1,\ldots,A_n$ be finite nonempty subsets of a field $F$, and
let
\begin{equation}\label{5.1}f(x_1,\ldots,x_n)=a_1x_1^k+\cdots+a_nx_n^k+g(x_1,\ldots,x_n)\in F[x_1,\ldots,x_n]\end{equation}
 with
 \begin{equation}\label{5.2}k\in\Z^+,\ a_1,\ldots,a_n\in F^*\ \t{and}\ \deg g<k.\end{equation}
 {\rm (i)} We have
 $$\begin{aligned}&|\{f(x_1,\ldots,x_n):\,x_1\in A_1,\ldots,x_n\in A_n\}|
 \\&\quad\gs\min\bg\{p(F),\,\sum_{i=1}^n\l\lfloor\f{|A_i|-1}k\r\rfloor+1\bg\}.
 \end{aligned}$$

 {\rm (ii)} If $k\gs n$ and $|A_i|\gs i$ for $i=1,\ldots,n$, then
 $$\begin{aligned} &|\{f(x_1,\ldots,x_n):\,x_1\in A_1,\ldots,x_n\in A_n,\ \t{and}\ x_i\not=x_j\ \t{if}\ i\not=j\}|
 \\&\qquad\qquad\gs\min\bg\{p(F),\,\sum_{i=1}^n\l\lfloor\f{|A_i|-i}k\r\rfloor+1\bg\}.
 \end{aligned}$$
\end{theorem}

\begin{remark} Let $a_1,\ldots,a_n$ be nonzero elements of a finite field $F$ and let $k$ be a positive integer.
Concerning lower bounds for
$|\{a_1x_1^k+\cdots+a_kx_n^k:\ x_1,\ldots,x_n\in F\}|$,
the reader may consult \cite{CMS} and \cite{T} for earlier results.
\end{remark}

Motivated by a concrete example, Sun \cite{S08a} actually raised the following extension of Conjecture \ref{Conj1.1}.

\begin{conjecture}[Sun \cite{S08a}]\label{Conj5.2} Let $f(x_1,\ldots,x_n)$
be a polynomial over a field $F$ given by $(5.1)$ and $(5.2)$.
Provided that $p(F)\not=n+1$ and $n>k$, for any finite subset $A$ of $F$ we
have
$$\begin{aligned}&|\{f(x_1,\ldots,x_n):\, x_1,\ldots,x_n\in A,\ \t{and}\ x_i\not=x_j\ \t{if}\ i\not=j\}|
\\\gs&\min\l\{p(F)-[\![n=2\ \&\ a_1=-a_2]\!],\, \f{n(|A|-n)-\{n\}_k\{|A|-n\}_k}k+1\r\},\end{aligned}$$
where we use $\{m\}_k$ to denote the least nonnegative residue of an integer $m$ modulo $k$.
\end{conjecture}

Sun \cite{S08a} proved the last inequality with
the lower bound replaced by $\min\{p(F),\,|A|-n+1\}$.

Theorem \ref{Th1.3} on restricted sumsets can be extended to the following general result on restricted value sets.

\begin{theorem}\label{Th5.2} Let $F$ be a field, and let
$f(x_1,\ldots,x_n)\in F[x_1,\ldots,x_n]$ be given by $(5.1)$ and $(5.2)$.
Let $P(x_1,\ldots,x_n)$ be a polynomial over $F$ with $[x_1^{k_1}\cdots x_n^{k_n}]P(x_1,\ldots,x_n)\not=0$,
where $k_1,\ldots,k_n$ are nonnegative integers with $k_1+\cdots+k_n=\deg P$.
Let $A_1,\ldots,A_n$ be finite subsets of $F$ with $|A_i|>k_i$ for $i=1,\ldots,n$.
Then, for the restricted value set
\begin{equation}\label{5.3}V=\{f(x_1,\ldots,x_n):\ x_1\in A_1,\ldots,x_n\in A_n,\ \t{and}\ P(x_1,\ldots,x_n)\not=0\},\end{equation}
we have
\begin{equation}\label{5.4}|V|\gs\min\bg\{p(F)-\sum_{i=1}^n\l\lfloor\f{k_i}k\r\rfloor,
\ \sum_{i=1}^n\(\l\lfloor\f{|A_i|-k_i-1}k\r\rfloor-\l\lfloor\f{k_i}k\r\rfloor\)+1\bg\}.\end{equation}
\end{theorem}
\Proof. It suffices to consider the nontrivial case
$$p(F)>\sum_{i=1}^n\l\lfloor\f{k_i}k\r\rfloor\ \t{and}\ \sum_{i=1}^n\(\l\lfloor\f{|A_i|-k_i-1}k\r\rfloor
-\l\lfloor\f{k_i}k\r\rfloor\)\gs0.$$

For $i=1,\ldots,n$ let $r_i$ be the least nonnegative residue of $k_i$ modulo $k$.
Write $P(x_1,\ldots,x_n)$ in the form (3.1) and
consider the polynomial
$$\bar P(x_1,\ldots,x_n)=\sum_{j_i\in r_i+k\N\ \t{for}\ i=1,\ldots,n
\atop j_1+\cdots+j_n=\deg P} c_{j_1,\ldots,j_n}\prod_{i=1}^na_i^{(r_i-j_i)/k}(x_i)_{(j_i-r_i)/k}.$$
Clearly
\begin{align*}&\bigg[\prod_{i=1}^nx_i^{\lfloor k_i/k\rfloor}\bigg]\bar P(x_1,\ldots,x_n)
\\=&\sum_{j_i\in k_i+k\N\ \t{for}\ i=1,\ldots,n
\atop \sum_{i=1}^nj_i=\sum_{i=1}^n k_i} c_{j_1,\ldots,j_n}\prod_{i=1}^na_i^{(r_i-j_i)/k}\cdot
\bigg[\prod_{i=1}^nx_i^{\lfloor k_i/k\rfloor}\bigg]\prod_{i=1}^n(x_i)_{(j_i-r_i)/k}
\\=&c_{k_1,\ldots,k_n}\prod_{i=1}^na_i^{(r_i-k_i)/k}\not=0.
\end{align*}
For $i=1,\ldots,n$ let $B_i=\{me:\ m\in I_i\}$ where
$$I_i=\bigg[\l\lfloor\f{|A_i|-r_i-1}k\r\rfloor-\l\lfloor\f{k_i}k\r\rfloor,
\ \l\lfloor\f{|A_i|-r_i-1}k\r\rfloor\bigg].$$
Clearly $|B_i|=\lfloor k_i/k\rfloor+1$ since $\lfloor k_i/k\rfloor<p(F)$.
Note also that
$$\l\lfloor\f{|A_i|-r_i-1}k\r\rfloor\gs\l\lfloor\f{k_i-r_i}k\r\rfloor=\f{k_i-r_i}k=\l\lfloor\f{k_i}k\r\rfloor.$$
In light of the Combinatorial Nullstellensatz, there are $q_1\in I_1,\ldots,q_n\in I_n$ such that
\begin{equation}\label{5.5}\bar P(q_1e,\ldots,q_ne)\not=0.\end{equation}

 Set $m_i=kq_i+r_i$ for $i=1,\ldots,n$. Then
 \begin{align*} M=&\f{\sum_{i=1}^n m_i-\deg P}k=\sum_{i=1}^n\f{m_i-k_i}k=\sum_{i=1}^n\(q_i-\l\lfloor\f{k_i}k\r\rfloor\)
 \\\gs&\sum_{i=1}^n\(\l\lfloor\f{|A_i|-k_i-1}k\r\rfloor-\l\lfloor\f{k_i}k\r\rfloor\)\gs0.
 \end{align*}
 and
 \begin{align*}&[x_1^{m_1}\cdots x_n^{m_n}]P(x_1,\ldots,x_n)f(x_1,\ldots,x_n)^M
 \\=&[x_1^{m_1}\cdots x_n^{m_n}]P(x_1,\ldots,x_n)(a_1x_1^k+\cdots+a_nx_n^k)^M
\\=&\sum_{j_i\in m_i-k\N\ \t{for}\ i=1,\ldots,n\atop j_1+\cdots+j_n=\deg P}
c_{j_1,\ldots,j_n}\f{M!}{\prod_{i=1}^n((m_i-j_i)/k)!}\prod_{i=1}^na_i^{(m_i-j_i)/k}
 \end{align*}
 So we have
 \begin{align*}&q_1!\cdots q_n![x_1^{m_1}\cdots x_n^{m_n}]P(x_1,\ldots,x_n)f(x_1,\ldots,x_n)^M
 \\=&M!\sum_{j_i\in m_i-k\N\ \t{for}\ i=1,\ldots,n\atop\ j_1+\cdots+j_n=\deg P}
c_{j_1,\ldots,j_n}\prod_{i=1}^na_i^{(m_i-j_i)/k}(q_ie)_{\lfloor j_i/k\rfloor}
\\=&M!a_1^{q_1}\cdots a_n^{q_n}\bar P(q_1e,\ldots,q_ne).
 \end{align*}

If $|V|\ls M<p(F)$, then by (\ref{5.5}) and the above we have
\begin{align*} &[x_1^{m_1}\cdots x_n^{m_n}]P(x_1,\ldots,x_n)f(x_1,\ldots,x_n)^{M-|V|}\prod_{v\in V}(f(x_1,\ldots,x_n)-v)
\\&\quad=[x_1^{m_1}\cdots x_n^{m_n}]P(x_1,\ldots,x_n)f(x_1,\ldots,x_n)^M\not=0,
\end{align*}
hence by the Combinatorial Nullstellensatz there are $x_1\in A_1,\ldots,x_n\in A_n$ such that
$$P(x_1,\ldots,x_n)f(x_1,\ldots,x_n)^{M-|V|}\prod_{v\in V}(f(x_1,\ldots,x_n)-v)\not=0$$
which contradicts (\ref{5.3}). Therefore, either
\begin{align*} p(F)\ls &M=\sum_{i=1}^n\(q_i-\l\lfloor\f{k_i}k\r\rfloor\)
\\\ls&
\sum_{i=1}^n\(\l\lfloor\f{|A_i|-r_i-1}k\r\rfloor-\l\lfloor\f{k_i}k\r\rfloor\)
=\sum_{i=1}^n\l\lfloor\f{|A_i|-k_i-1}k\r\rfloor
\end{align*}
 or $$|V|\gs M+1\gs \sum_{i=1}^n\(\l\lfloor\f{|A_i|-k_i-1}k\r\rfloor-\l\lfloor\f{k_i}k\r\rfloor\)+1.$$

If $p(F)>\sum_{i=1}^n\lfloor(|A_i|-k_i-1)/k\rfloor$, then we have
\begin{align*}|V|\gs&\sum_{i=1}^n\(\l\lfloor\f{|A_i|-k_i-1}k\r\rfloor-\l\lfloor\f{k_i}k\r\rfloor\)+1
\\=&\min\bg\{p(F)-\sum_{i=1}^n\l\lfloor\f{k_i}k\r\rfloor,
\ \sum_{i=1}^n\(\l\lfloor\f{|A_i|-k_i-1}k\r\rfloor-\l\lfloor\f{k_i}k\r\rfloor\)+1\bg\}.
\end{align*}

 In the case $p(F)\ls \sum_{i=1}^n\lfloor(|A_i|-k_i-1)/k\rfloor$, as $\sum_{i=1}^n k_i=\deg P$
there are $A_1'\se A_1,\ldots,A_n'\se A_n$ such that
$$|A_1'|>k_1,\ldots,|A_n'|>k_n,\ \t{and}\ \sum_{i=1}^n\l\lfloor\f{|A_i'|-k_i-1}k\r\rfloor=p(F)-1<p(F),$$
therefore
\begin{align*} |V|\gs&|\{x_1+\cdots+x_n:\ x_1\in A_1',\ldots,x_n\in A_n',\ \t{and}\ P(x_1,\ldots,x_n)\not=0\}|
\\\gs&\min\bg\{p(F)-\sum_{i=1}^n\l\lfloor\f{k_i}k\r\rfloor,\
\sum_{i=1}^n\(\l\lfloor\f{|A_i'|-k_i-1}k\r\rfloor-\l\lfloor\f{k_i}k\r\rfloor\)+1\bg\}
\\=&p(F)-\sum_{i=1}^n\l\lfloor\f{k_i}k\r\rfloor
\\=&\min\bg\{p(F)-\sum_{i=1}^n\l\lfloor\f{k_i}k\r\rfloor,
\ \sum_{i=1}^n\(\l\lfloor\f{|A_i|-k_i-1}k\r\rfloor-\l\lfloor\f{k_i}k\r\rfloor\)+1\bg\}.
\end{align*}
We are done. \qed

 Here is a consequence of Theorem \ref{Th5.2}.

\begin{corollary}\label{Cor5.1} Let $F$ be a field and let $f(x_1,\ldots,x_n)\in F[x_1,\ldots,x_n]$ be
given by $(\ref{5.1})$ and $(\ref{5.2})$.
Let
$A_1,\ldots,A_n$ be finite subsets of $F$ with $|A_i|\gs i$ for $i=1,\ldots,n$.
Then, for the restricted value set
\begin{equation}\label{5.6} V=\{f(x_1,\ldots,x_n):\ x_1\in A_1,\ldots,x_n\in A_n,\ \t{and}\ x_1,\ldots,x_n\ \t{are distinct}\},
\end{equation}
we have
\begin{equation}\label{5.7}|V|+\Delta(n,k)\gs\min\bg\{p(F),
\sum_{i=1}^n\l\lfloor\f{|A_i|-i}k\r\rfloor+1\bg\},\end{equation}
where
\begin{equation}\label{5.8}\Delta(n,k)=\l\lfloor\f nk\r\rfloor\(n-k\f{\lfloor n/k\rfloor+1}2\).\end{equation}
\end{corollary}
\Proof. We apply Theorem \ref{Th5.2} with
\begin{align*} P(x_1,\ldots,x_n)=&\prod_{1\ls i<j\ls n}(x_j-x_i)=\det(x_j^{i-1})_{1\ls i,j\ls n}.
\end{align*}
Note that $[\prod_{i=1}^nx_i^{i-1}]P(x_1,\ldots,x_n)=1\not=0$.
By Theorem \ref{Th5.2},
$$|V|+\sum_{i=1}^n\l\lfloor\f{i-1}k\r\rfloor\gs\min\bg\{p(F),
\ \sum_{i=1}^n\l\lfloor\f{|A_i|-i}k\r\rfloor+1\bg\}.$$
So it suffices to observe that
\begin{align*} \sum_{i=1}^n\l\lfloor\f{i-1}k\r\rfloor
=&\sum_{q=0}^{\lfloor n/k\rfloor-1}\sum_{r=1}^k\l\lfloor\f{qk+r-1}k\r\rfloor
+\sum_{k\lfloor n/k\rfloor<i\ls n}\l\lfloor\f{i-1}k\r\rfloor
\\=&\sum_{q=0}^{\lfloor n/k\rfloor-1}kq+\(n-k\l\lfloor\f nk\r\rfloor\)\l\lfloor\f nk\r\rfloor
\\=&k\l\lfloor\f nk\r\rfloor\f{\lfloor n/k\rfloor-1}2+\(n-k\l\lfloor\f nk\r\rfloor\)\l\lfloor\f nk\r\rfloor
=\Delta(n,k).
\end{align*}
This concludes the proof. \qed

\begin{lemma}\label{Lem5.1} Let $k$ and $n$ be positive integers. Then, for any $m\in\Z$ we have
\begin{equation}\label{5.9}\begin{aligned}\sum_{i=1}^n\l\lfloor\f{m-i}k\r\rfloor=&m\l\lfloor\f nk\r\rfloor+\{n\}_k\l\lfloor\f{m-n}k\r\rfloor
-\f k2\l\lfloor\f nk\r\rfloor\(\l\lfloor\f nk\r\rfloor+1\)
\\&+\{m\}_k[\![\{m\}_k<\{n\}_k]\!].
\end{aligned}\end{equation}
\end{lemma}
\Proof. Let $f(m)$ and $g(m)$ denote the left-hand side and the right-hand side of (\ref{5.9}) respectively.
We first prove that $f(n)=g(n)$. In fact, by the proof of Corollary \ref{Cor5.1},
$$f(n)=\sum_{j=0}^{n-1}\l\lfloor\f{j}n\r\rfloor=\Delta(n,k)=g(n).$$

Next we show that $f(m+1)-f(m)=g(m+1)-g(m)$ for any $m\in\Z$. Observe that
\begin{align*} f(m+1)-f(m)=&\sum_{i=1}^n\(\l\lfloor\f{m+1-i}k\r\rfloor-\l\lfloor\f{m-i}k\r\rfloor\)
\\=&|\{1\ls i\ls n:\ i\eq m+1\ (\mo\ k)\}|
\\=&|\{q\in\N:\ \{m\}_k+kq<n\}|=\l\lfloor\f nk\r\rfloor+[\![\{m\}_k<\{n\}_k]\!].
\end{align*}
Also,
\begin{align*} &g(m+1)-g(m)-\l\lfloor\f nk\r\rfloor
\\=&\{n\}_k[\![m+1\eq n\ (\mo\ k)]\!]
\\&+\{m+1\}_k[\![\{m+1\}_k<\{n\}_k]\!]-\{m\}_k[\![\{m\}_k<\{n\}_k]\!]
\\=&\{m+1\}_k[\![\{m+1\}_k\ls\{n\}_k]\!]-\{m\}_k[\![\{m\}_k<\{n\}_k]\!]
=[\![\{m\}_k<\{n\}_k]\!].
\end{align*}

 So far we have proved (\ref{5.9}) for all $m\in\Z$. \qed

The following result partially resolves Conjecture 5.1.

\begin{corollary}\label{Cor5.2} Let $F$ be a field and let $f(x_1,\ldots,x_n)\in F[x_1,\ldots,x_n]$
be given by $(\ref{5.1})$ and $(\ref{5.2})$.
Let $A_1,\ldots,A_n$ be finite subsets of $F$ with $|A_1|=\cdots=|A_n|=m\gs n$.
Then, for the restricted value set $V$ in $(\ref{5.6})$
we have
\begin{equation}\label{5.10}|V|\gs\min\bg\{p(F)-\Delta(n,k),
\f{n(m-n)-\{n\}_k\{m-n\}_k}k+r_{k,m,n}+1\bg\},\end{equation}
where
\begin{equation}\label{5.11}r_{k,m,n}=\{m\}_k[\![\{m\}_k<\{n\}_k]\!].\end{equation}
\end{corollary}

\begin{remark} In the special case $a_1=\cdots=a_n$, H. Pan and Sun \cite{PS2} proved (\ref{5.10}) with $\Delta(n,k)$ omitted.
\end{remark}

\noindent{\it Proof of Corollary \ref{Cor5.2}}. By Lemma \ref{Lem5.1},
\begin{align*} &\sum_{i=1}^n\l\lfloor\f{m-i}k\r\rfloor-\Delta(n,k)
\\=&(m-n)\l\lfloor\f nk\r\rfloor+\{n\}_k\l\lfloor\f{m-n}k\r\rfloor+r_{k,m,n}
\\=&\f{n(m-n)}k-\{n\}_k\f{m-n}k+\{n\}_k\l\lfloor\f{m-n}k\r\rfloor+r_{k,m,n}
\\=&\f{n(m-n)-\{n\}_k\{m-n\}_k}k+r_{k,m,n}.
\end{align*}
So, the desired result follows from Corollary \ref{Cor5.1}. \qed

\Ack. The authors are grateful to the two referees for their helpful comments.

\end{document}